\documentclass{article}
\usepackage{amssymb}
\usepackage{graphicx}
\usepackage{amsfonts}
\usepackage{latexsym}
\usepackage{amsthm}
\usepackage{amsmath}

\usepackage{amssymb, amscd}

\newtheorem{problem}{Problem}[section]
\newtheorem{theo}[problem]{Theorem}
\newtheorem{rem}[problem]{Remark}

\newtheorem{defin}[problem]{Definition}
\newtheorem{prop}[problem]{Proposition}
\newtheorem{cor}[problem]{Corollary}

\begin{document}
\date{August 2004}
\title{{\Large Arrangements of symmetric products of spaces}}

\author{Pavle Blagojevi\' c \and Vladimir Gruji\' c \and  \and Rade \v
Zivaljevi\' c\footnote{The authors were supported by the grant
1643 of the Serbian Ministry of Science, Technology and
Development.}}

\maketitle

\begin{abstract}
We study the combinatorics and topology of general arrangements of
subspaces of the form $D + SP^{n-d}(X)$ in symmetric products
$SP^n(X)$ where $D\in SP^d(X)$. Symmetric products $SP^m(X):=
X^m/S_m$, also known as the spaces of  effective ``divisors'' of
order $m$, together with their companion spaces of
divisors/particles, have been studied from many points of view in
numerous papers, see \cite{BGZ} and \cite{Kallel_1998} for the
references. In this paper we approach them from the point of view
of geometric combinatorics. Using the topological technique of
{\em diagrams of spaces} along the lines of \cite{WZZ99} and
\cite{ZieZiv93}, we calculate the homology of the union and the
complement of these arrangements. As an application we include a
computation of the homology of the {\em homotopy end space} of the
open manifold $SP^n(M_{g,k})$, where $M_{g,k}$ is a Riemann
surface of genus $g$ punctured at $k$ points, a problem which was
originally motivated by the study of commutative $(m+k,m)$-groups
\cite{TD2001}.
\end{abstract}

\section{Arrangements of symmetric products}
\label{sec:arrangement}

The study of homotopy types of arrangements of subspaces with an
emphasis on the underlying combinatorial structure is a well
establish part of geometric and topological combinatorics.
Originally the focus was on the arrangements of linear or affine
subspaces \cite{BjoGu} \cite{Gru} \cite{OrlTer92} \cite{Zas75}.
 Gradually other, more general arrangements of spaces were
introduced and studied. Examples include arrangements of
pseudolines/pseudospheres, in connection with the realizations of
(oriented) matroids \cite{BLVSWZ} \cite{Swartz}, arrangements of
projective and Grassmann varieties, partially motivated by a
geometrization of the Stanley ring construction \cite{GoMacPh}
\cite{WZZ99}, arrangements of classifying spaces $BH$ for a family
of subgroups of a given group, subspace arrangements over finite
fields \cite{BjoEke97} etc. With the introduction into
combinatorics of the technique of {\em diagrams of spaces} and the
associated {\em homotopy colimits}  \cite{WZZ99} \cite{ZieZiv93},
it became apparent that arrangements of subspaces have much in
common with other important and well studied objects like
stratified spaces/discriminants and their geometric resolutions
\cite{Vas94}, toric varieties viewed as combinatorial objects
associated to face latices of polytopes \cite{BuPan99}
\cite{DaJan91} \cite{WZZ99} etc. All this serves as a motivation
for the study of general subspace arrangements carrying
interesting combinatorial structure.

\medskip
Symmetric products of spaces $SP^n(X)$ are classical mathematical
objects  \cite{ACGH} \cite{Arnold_96} \cite{C2M2} \cite{Dijk}
\cite{DiMoVeVe} \cite{DoldThom} \cite{Kallel_1998} \cite{Segal}
which appear in different areas of mathematics and mathematical
physics as orbit spaces, divisor spaces, particle spaces etc., see
\cite{BGZ} for a leisurely introduction and a review of old and
new applications. The case of $2$-manifolds $M$ is of particular
interest since in this case $SP^n(M)$ is a manifold. Elements
$D\in SP^d(M)$ are called divisors of order $\vert D\vert = d$. In
this paper we study arrangements in $SP^n(X)$ of the form
\begin{equation}
\label{eqn:aranzman}
 {\mathcal A} = \{D_i + SP^{n-\vert
D_i\vert}(X)\}_{i=1}^k ,
\end{equation}
where the case of open or closed surfaces is of special interest.

\medskip
Given an arrangement ${\mathcal A} = \{F_1,\ldots, F_k\}$ of
subspaces in an ambient space $V$, the {\em union} or the {\em
link} of ${\mathcal A}$ is $D({\mathcal A}) :=
\bigcup_{i=1}^k~F_i$ and the complement is $M({\mathcal A}) :=
V\setminus D({\mathcal A})$. In this paper we compute the homology
of the link and, in the case of Riemann surfaces $M_g$, the
homology of the complement of the arrangement
(\ref{eqn:aranzman}). As an application we compute the homology of
the ``homotopy-end-space'' of $SP^n(M_{g,k})$ where $M_{g,k} :=
M_g\setminus \{x_1,\ldots ,x_k\}$ is the so called $(g,k)$-amoeba,
see Figure~\ref{fig:amoeba}, and discuss the connection of this
computation with the problem of existence of commutative
$(m+k,m)$-groups and the problem what information about the
original surface $M$ can be reconstructed from the symmetric
product $SP^n(M)$, Section~\ref{sec:commutative}.

Note that if $g=0$ and if $D_i = x_i$ are all distinct divisors of
order $1$, then the complement of the arrangement
(\ref{eqn:aranzman}) is homeomorphic to the complement of a
generic arrangement of $k-1$ hyperplanes in $\mathbb{C}^n$,
studied by Hattori \cite{Hat75} \cite{OrlTer92}. So the results
about $SP^n(M_{g,k})$ can be viewed as an extension of some
classical results about complex hyperplane arrangements. On the
technical side, we would like to emphasize the role of the
relation of {\em proper} domination between simple diagrams,
Definition~\ref{def:dominate} in Section~\ref{sec:homol}, which
often allows, as demonstrated in Section~\ref{sec:category}, the
study of naturality properties for Goresky-MacPherson
\cite{GoMacPh} and Ziegler--{\v Z}ivaljevi{\' c} \cite{ZieZiv93}
\cite{WZZ99} type formulas. The same concept allows us to extend
(Theorems~\ref{thm:suma} and \ref{thm:opetsuma}) the classical
Steenrod's theorem on the decompositions od symmetric products, to
the case od diagrams of spaces.

\subsection{Homology of the union $D(\mathcal{A}) = \bigcup \mathcal{A}$}
\label{sec:homol}

We approach the computation of the homology
$H_{\ast}(D(\mathcal{A}))$ by the method of {\em diagrams of
spaces}. The references emphasizing applications of this technique
in geometric combinatorics are \cite{WZZ99} and \cite{ZieZiv93}.
The reader is referred to these papers for the notation and
standard facts. Aside from standard tools like the {\em Projection
Lemma} or the {\em Homotopy Lemma}, we make a special use of the
idea of an ``ample space'' diagram outlined in Section~5.4., of
\cite{WZZ99}.

\medskip
Let ${\mathcal A}=\{F_n^i\}_{i=1}^k$ be an arrangement of
subspaces in $SP^n({X})$ where $F_n^i := x_i + SP^{n-1}({X})\cong
SP^{n-1}({X})$ and $\{x_i\}_{i=1}^k$ is a collection of distinct
points in ${X}$. Let $P = P({\mathcal A})$ be the intersection
poset of ${\mathcal A}$. By definition \cite{OrlTer92}
 \cite{WZZ99} \cite{ZieZiv93}, $P$ has an element for each
non-empty intersection $F_n^I = F_n^{i_0}\cap
F_n^{i_1}\cap\ldots\cap F_n^{i_p}$ where $I = \{i_0,i_1,\ldots,
i_p\}\subset [k]:=\{1,\ldots, k\}.$ If $n\geq k$ then $P$ is
isomorphic to the power set $\mathcal{P}'[k] =
\mathcal{P}[k]\setminus\{\emptyset\}$ or alternatively, the face
poset of an abstract simplex $\Sigma$ with vertices $\{1,\ldots,
k\}$. If $n\leq k$ then $P$ is isomorphic to the poset
$\mathcal{P}'_{\leq n}[k]$ of all non-empty subsets $I$ of
cardinality at most $n$ or alternatively the $(n-1)$-skeleton
$\Sigma^{n-1}$ of $\Sigma$.

More generally assume that $\mathcal{A}=\{F_n^i\}_{i=1}^r$ is an
arrangement of subspaces of the form
\begin{equation}
\label{eqn:general}
 F_n^i = D_i + SP^{n-d_i}({X}) \quad \mbox{ \rm where } \quad
 D_i\in SP^{d_i}({X}) \quad i=1,\ldots, r.
\end{equation}
Moreover, we assume that for each $i$ the corresponding divisor
$D_i$ has the form
\[
D_i = \alpha^i_1x_1+\ldots + \alpha^i_kx_k
\]
where all points $x_i\in {X}$ are distinct and fixed in advance
while $\alpha_i$ are non-negative integers. The associated
intersection poset has several useful interpretations. Let $J =
\{j_1<\ldots <j_m\}$ be a subsequence of $[k] = \{1,\ldots, k\}$.
Since $K\in \bigcap_{\alpha=1}^m~(D_{j_\alpha} +
SP^{n-d_{j_\alpha}}({X}))$ if and only if $K\in SP^n({X})$ and
$D_{j_\alpha}\leq K$ for each $\alpha =1,\ldots,m$ we observe that
$K\geq D_{j_1,j_2,\ldots, j_m} =: D_J$ where $D_J$ is the least
upper bound of divisors $D_{j_\alpha}, \, \alpha=1,\ldots, m$.
Another description is in terms of {\em multisets} \cite{Stanley}
or the associated monomials. Given an effective divisor $D =
\alpha_1x_1+\ldots +\alpha_rx_r$, where points $x_i$ are distinct
and $\alpha_i$ are the corresponding multiplicities, the
associated multiset (monomial) is $x_1^{\alpha_1}\ldots
x_r^{\alpha_r}$. Then the intersection poset $P = P(\mathcal{A})$
can be described as the collection of all multisets in ${X}$ of
cardinality at most $n$ which can be represented as unions of
multisets associated to original divisors $D_i$. Yet another
description arises if multisets $x_1^{\alpha_1}\ldots
x_r^{\alpha^r}$ are interpreted as natural numbers
$p_1^{\alpha_1}\ldots p_r^{\alpha_r}$ where $p_i$ are distinct
prime integers. This shows that the combinatorics of intersections
posets of arrangements of subspaces in $SP^n(X)$ is directly
connected with the classical elementary number theory.

\medskip\noindent
{\bf Caveat:} All spaces we deal with are {\em admissible} in the
sense of Definition~\ref{def:admissible}. Elements of the
intersection poset $P = P(\mathcal{A})$ are often denoted by
$p,q,r$ etc.\ but when we want to emphasize that they are actually
divisors (multisets) we use the notation $I, J, K, D_j$ etc.

\medskip
Let $\mathcal{A}$ be an arrangement of subspaces described by
(\ref{eqn:general}). Let $\mathcal{D} : P \rightarrow Top$ be an
associated diagram of spaces and inclusion maps, \cite{WZZ99}
\cite{ZieZiv93}. Each $I\in P$ is of the form $I =
\beta_1x_1+\ldots +\beta_kx_k$ where $\vert I\vert =
\beta_1+\ldots +\beta_k\leq n$. Hence,
\[
\mathcal{D}(I) := \beta_1x_{1}+\ldots + \beta_kx_{k} + SP^{n-\vert
I\vert}({X}).
\]
By Projection Lemma, \cite{WZZ99} Lemma~4.5, or by
Proposition~6.9.\ on page 49 in \cite{Bal},
\begin{equation}
\label{eqn:projection} F_n = \bigcup
\mathcal{A}\cong\mathbf{colim}~\mathcal{D}  \simeq
\mathbf{hocolim}~\mathcal{D} .
\end{equation}

\begin{defin}
\label{def:sum} For $A\subset SP^p(X)$ and $B\subset SP^q(X)$, the
``Minkowski'' sum $A+B$ is a subset of $SP^{p+q}(X)$ defined by
$A+B := \{a+b \mid a\in A, b\in B\}.$
\end{defin}
Since $X$ is an admissible space, Definition~\ref{def:admissible},
there exists a closed, contractible set $C\supset\{x_1,\ldots,
x_k\}$ such that the projection map $X\rightarrow X/C$ is a
homotopy equivalence. Moreover, $C$ contracts to a point $y\in C$
which can be prescribed in advance.

 \medskip
Define $\mathcal{E} : P \rightarrow Top$ to be the diagram of
spaces and inclusion maps determined by
\[
\mathcal{E}(I) := \vert I\vert y + SP^{n-\vert I\vert}({X}) \quad
\mbox{ for each } \quad I \in P.
\]
Note that $\mathcal{E}(I)$ depends only on the order $\vert
I\vert$ of the divisor $I$. We would like to show that
\[
\mathbf{hocolim}~\mathcal{D} \simeq \mathbf{hocolim}~\mathcal{E} .
\]
Since there does not exist an obvious map between these diagrams,
we define a new diagram $\mathcal{C}$, a so called ``ample space''
diagram, which contains both $\mathcal{D}$ and $\mathcal{E}$ as
subdiagrams.

Let $\mathcal{C} : P \rightarrow Top$ be the diagram of spaces and
inclusion maps defined by
\[
\mathcal{C}(I) := SP^{\vert I\vert}(C) + SP^{n-\vert I\vert}({X})
\quad \mbox{ for each } \quad I  \in P.
\]

\begin{prop}
\label{prop:ample} Let
\[
\alpha : \mathcal{D}\rightarrow \mathcal{C} \quad \mbox{ and
}\quad \beta : \mathcal{E}\rightarrow\mathcal{C}
\]
be the morphisms of diagrams where $\alpha_I :
\mathcal{D}(I)\rightarrow \mathcal{C}(I)$ and $\beta_I :
\mathcal{E}(I)\rightarrow \mathcal{C}(I)$ are obvious inclusions.
Then $\alpha$ and $\beta$ induce the homotopy equivalences of the
corresponding homotopy colimits,
\[
\mathbf{hocolim}~\mathcal{D}\stackrel{\hat\alpha}{\longrightarrow}
\mathbf{hocolim}~\mathcal{C}\stackrel{\hat\beta}{\longleftarrow}
\mathbf{hocolim}~\mathcal{E} .
\]
\end{prop}

 \noindent
{\bf Proof:} The proposition is an immediate consequence of the
Homotopy Lemma, \cite{WZZ99} Lemma~4.6, and
Proposition~\ref{prop:main} from Section~\ref{sec:main}.

\subsection{Steenrod's theorem for diagrams}\label{sec:steenrod}

In the previous section the calculation of the homology of the
union of the arrangement ${\mathcal A}=
\{D_i+SP^{n-d_i}(X)\}_{i=1}^r$ was reduced to the calculation of
$H_\ast(\mathbf{hocolim}~\mathcal{E})$ for a diagram $\mathcal{E}$
of particularly simple form, cf.\ Definition~\ref{def:simple}.

Our objective in this section is to establish a decomposition
result (Theorem~\ref{thm:suma}) expressing the homology of
$\mathbf{hocolim}~\mathcal{E}$ in simpler terms.
Theorem~\ref{thm:suma} can be seen as an offspring and a
generalization of the well known Steenrod's theorem,
\[
H_\ast(SP^m(X);\mathbb{A}) \cong \bigoplus_{j=0}^m H_\ast(SP^j(X),
SP^{j-1}; \mathbb{A})
\]
where by definition $SP^{-1}(X)=\emptyset$ and $\mathbb{A}$ is an
arbitrary Abelian group.
\begin{defin}
\label{def:standard} Assume that $y\in {X}$ is a point fixed in
advance. A {\em standard} inclusion \[e_{p,q} :
SP^p({X})\hookrightarrow SP^q({X}),\] associated to a pair of
integers $p\leq q$, is by definition the map defined by
$e_{p,q}(Y) := (q-p)y + Y$ for each $Y\in SP^p({X})$.
\end{defin}

\begin{defin}
\label{def:simple} A diagram of spaces $\mathcal{D} : P\rightarrow
Top$ is called {\em simple} if
\begin{itemize}
 \item[{\rm (a)}]  for each $p\in P$ there exists $\mu=\mu(p)$ such that $\mathcal{D}(p) =
 SP^{\mu}({X})$,
 \item[{\rm (b)}] the map $\mathcal{D}_{q,p} :
SP^{\mu(p)}(X) \rightarrow SP^{\mu(q)}(X)$ is a standard inclusion
for each pair $q\leq p$.
\end{itemize}
The monotone function $\mu : P\rightarrow \mathbb{N}$ is called
the {\em rank} function of $\mathcal{D}$. If $\mu(p) = \mu(q)$ for
each pair $p,q\in P$ we say that $\mathcal{D}$ is a {\em constant}
diagram.
\end{defin}

{\rm
\begin{rem}
\label{rem:simple} Note that a diagram of spaces $\mathcal{D} :
P\rightarrow Top$ is simple if and only if there exists a strictly
increasing sequence $m = c_0 < c_1 <\ldots <c_k = M$ and a
strictly decreasing sequence $P_m = P_{c_0} \supset
P_{c_1}\supset\ldots\supset P_{c_k}= P_M$ of ideals in $P$ such
that ${\rm Im}(\mu) = \{c_i\}_{i=0}^k$ and for each $p\in
P_{c_i}\setminus P_{c_{i+1}}, \, \mathcal{D}(p) = SP^{c_i}({X})$.
\end{rem}
}

\begin{defin}
\label{def:dominate} Suppose that $\mathcal{E}_1$ and
$\mathcal{E}_2$ are both simple diagrams over the same poset $P$
and let $\mu_1$ and $\mu_2$ be the corresponding rank functions.
We write $\mathcal{E}_1\preccurlyeq\mathcal{E}_2$ and say that the
diagram $\mathcal{E}_1$ is dominated by the diagram
$\mathcal{E}_2$ if $\mu_1\leq\mu_2$. We write
$\mathcal{E}_1\preccurlyeq_P\mathcal{E}_2$ and say that the
diagram $\mathcal{E}_1$ is {\em properly} dominated by
$\mathcal{E}_2$ if $\mu_1 = \mbox{\rm min}\{\mu_2, c\}$ for some
constant $c\in \mathbb{N}$.
\end{defin}

It is obvious that both $\preccurlyeq$ and $\preccurlyeq_P$ are
partial orders on the set of all simple diagrams over $P$. If
either $\mathcal{E}_1\preccurlyeq\mathcal{E}_2$ or
$\mathcal{E}_1\preccurlyeq_P\mathcal{E}_2$, there is a unique
morphism $\alpha : \mathcal{E}_1\rightarrow\mathcal{E}_2$ such
that $\alpha_p : \mathcal{E}_1(p)\rightarrow \mathcal{E}_2(p)$ is
a standard inclusion for each $p\in P$.

\begin{prop}\label{prop:injective} Suppose that
$\mathcal{E}_1\preccurlyeq_P\mathcal{E}_2\preccurlyeq_P\mathcal{E}_3$
and let $\mathcal{E}_1\stackrel{\alpha}{\longrightarrow}
 \mathcal{E}_2\stackrel{\beta}{\longrightarrow}
 \mathcal{E}_3$ be the associated chain of morphisms. Suppose that
 \[
 X_{\mathcal{E}_1}\stackrel{\hat{\alpha}}{\longrightarrow}
 X_{\mathcal{E}_2}\stackrel{\hat{\beta}}{\longrightarrow}
 X_{\mathcal{E}_3}
 \]
 is the corresponding chain of homotopy
 colimits $X_{\mathcal{E}_i}:= {\bf hocolim}~\mathcal{E}_i, \,
 i=1,2,3$. Then the map
 \[
 H_\ast(X_{\mathcal{E}_2}, X_{\mathcal{E}_1}; \mathbb{A})
 \longrightarrow H_\ast(X_{\mathcal{E}_3}, X_{\mathcal{E}_1}; \mathbb{A})
 \]
 is injective and the associated long exact sequence of the triple splits.
\end{prop}

\medskip\noindent
{\bf Proof:} It is sufficient to prove that the map
\begin{equation}
\label{eqn:inject} \hat\alpha : H_\ast(X_\mathcal{D};
\mathbb{A})\longrightarrow H_\ast(X_\mathcal{E}; \mathbb{A})
\end{equation}
is injective for each pair of simple diagrams such that
$\mathcal{D}\preccurlyeq_P\mathcal{E}$. Indeed, on applying this
result on pairs $\mathcal{E}_1\preccurlyeq_P\mathcal{E}_2$ and
$\mathcal{E}_1\preccurlyeq_P\mathcal{E}_3$ we obtain a commutative
diagram with exact rows,

 {\small
\[
\begin{CD}
0 @>>> H_\ast(X_{\mathcal{E}_1}) @>>> H_\ast(X_{\mathcal{E}_2})
@>>> H_\ast(X_{\mathcal{E}_2}, X_{\mathcal{E}_1}) @>>> 0 \\
 @. @VVV @VVV @VVV\\
 0 @>>> H_\ast(X_{\mathcal{E}_1}) @>>> H_\ast(X_{\mathcal{E}_3})
@>>> H_\ast(X_{\mathcal{E}_3}, X_{\mathcal{E}_1}) @>>> 0
\end{CD}
\]
 }

The first vertical arrow from the left is an isomorphism and the
second is a monomorphism, so the result follows from a version of
$5$-lemma.

\medskip
Let us observe that (\ref{eqn:inject}) is obvious if both
$\mathcal{D}$ and $\mathcal{E}$ are constant simple diagrams in
the sense that for some integers $m\leq n$, $\mathcal{D}(p) =
SP^m({X})$ and $\mathcal{E}(p) = SP^n({X})$ for each $p\in P$.
Indeed, in this case (\ref{eqn:inject}) reduces to the
monomorphism

 {\small
\begin{equation}
\label{eqn:mono}
 H_\ast(SP^m({X})\times \Delta(P))\longrightarrow
 H_\ast(SP^n({X})\times \Delta(P))) .
\end{equation}
 }
In light of the fact that by Steenrod's theorem
$H_\ast(SP^m({X}))\rightarrow H_\ast(SP^n({X}))$ is always a
monomorphism to a direct summand of $H_\ast(SP^n({X}))$,
(\ref{eqn:mono}) follows from K\" unneth formula and the
$5$-lemma. Next we observe that (\ref{eqn:inject}) is true even if
only $\mathcal{D}$ is a constant diagram. Indeed, let
$\mathcal{F}$ be a constant simple diagram such that
$\mathcal{E}\preccurlyeq\mathcal{F}$. Then the composition
$\hat{\beta}\circ\hat{\alpha}$ in the diagram
\[
\begin{CD}
 H_\ast(X_\mathcal{D}) @>\hat{\alpha}>> H_\ast(X_\mathcal{E})
 @>\hat{\beta}>> H_\ast(X_\mathcal{F})
\end{CD}
\]
is a monomorphism, hence $\hat\alpha$ alone is also a
monomorphism.

The general case of (\ref{eqn:inject}) is established by induction
on the size of the poset $P$. Let $\mathcal{C}$ be a maximal
constant simple diagram over $P$ such that
$\mathcal{C}\preccurlyeq\mathcal{D}$. In other words if
$\mathcal{D}(p) = SP^{\mu(p)}({X})$ for each $p\in P$ then
$\mathcal{C}(p):= SP^m({X})$ where $m:=\mbox{min}\{\mu(p)\}_{p\in
P}$. Let $P'$ be a subposet of $P$ defined by $P' := \{p\in P \mid
\mu(p)>m\}$. Note that $P'$ is actually an ideal in $P$. Define
$\mathcal{C}', \mathcal{D}',\mathcal{E}'$ respectively as the
restrictions of diagrams $\mathcal{C}, \mathcal{D}, \mathcal{E}$
on the subposet $P'$. Then by the excision axiom there is a
commutative diagram of long exact sequences
 {\small
\begin{equation}
\label{eqn:excision}
\begin{CD}
\ldots @>>> H_\ast(X_{\mathcal{D}'}, X_{\mathcal{C}'}) @>>>
H_\ast(X_{\mathcal{E}'}, X_{\mathcal{C}'}) @>>>
H_\ast(X_{\mathcal{E}'}, X_{\mathcal{D}'}) @>>> \ldots\\
@. @V\cong VV @V\cong VV @V\cong VV @.\\
\ldots @>>> H_\ast(X_{\mathcal{D}}, X_{\mathcal{C}}) @>>>
H_\ast(X_{\mathcal{E}}, X_{\mathcal{C}}) @>>>
H_\ast(X_{\mathcal{E}}, X_{\mathcal{D}}) @>>> \ldots
\end{CD}
\end{equation}
 }

The condition
$\mathcal{C}\preccurlyeq_P\mathcal{D}\preccurlyeq_P\mathcal{E}$
implies
$\mathcal{C}'\preccurlyeq_P\mathcal{D}'\preccurlyeq_P\mathcal{E}'$
and by the inductive assumption the first row splits. Hence there
is a short exact sequence
 {\small
\[
\begin{CD}
0 @>>> H_\ast(X_\mathcal{D},X_\mathcal{C}) @>\mbox{1--1} >>
H_\ast(X_\mathcal{E},X_\mathcal{C}) @>>>
H_\ast(X_\mathcal{E},X_\mathcal{D}) @>>> 0 .
\end{CD}
\]
 }
Finally, from the commutative diagram
 {\small
\[
 \begin{CD}
0 @>>> H_\ast(X_\mathcal{C}) @>\cong >> H_\ast(X_\mathcal{C}) @>>>
0\\
@. @VVV @VVV @VVV @.\\
0 @>>> H_\ast(X_\mathcal{D}) @ >\hat{\alpha} >>
H_\ast(X_\mathcal{E}) @>>>
H_\ast(X_\mathcal{E}, X_\mathcal{D}) @>>> 0\\
@. @VVV @VVV @VVV @.\\
0 @>>> H_\ast(X_\mathcal{D}, X_\mathcal{C}) @ >\mbox{1--1}>>
H_\ast(X_\mathcal{E},X_\mathcal{C}) @>>> H_\ast(X_\mathcal{E},
X_\mathcal{D}) @>>> 0
\end{CD}
\]
 }
we deduce that $\hat{\alpha}$ is a monomorphism. \hfill $\square$

\medskip

Given a simple diagram $\mathcal{D} : P\rightarrow Top, \,
p\mapsto SP^{\mu(p)}({X})$, let $m := \mbox{min}\{\mu(p)\}_{p\in
P}$ and $M := \mbox{max}\{\mu(p)\}_{p\in P}$. Let us assemble the
elements of the set $\{\mu(p)\}_{p\in P}$ into an increasing
sequence
\begin{equation}
\label{eqn:assemble} m = c_0 < c_1 <\ldots <c_k = M .
\end{equation}
Define $P_j := \{p\in P \mid \mu(p)\geq j\}$. Note that
\begin{equation}
\label{eqn:ideals} P = P_m \supseteq
P_{m+1}\supseteq\ldots\supseteq P_M
\end{equation}
is a decreasing sequence of ideals in $P$. In light of
(\ref{eqn:assemble}) we observe that
\begin{equation}
\label{egn:strict} P_m = P_{c_0} \supset
P_{c_1}\supset\ldots\supset P_{c_k}= P_M
\end{equation}
is a subsequence of (\ref{eqn:ideals}) obtained by removing the
redundant posets.

\begin{theo}
\label{thm:suma} Assume that $\mathcal{D} : P\rightarrow Top$ is a
simple diagram where $\mu : P \rightarrow \mathbb{N}$ is the
corresponding rank function, $m= \mbox{\rm min}\{\mu(p)\}_{p\in
P}, M= \mbox{\rm max}\{\mu(p)\}_{p\in P}$ and $P_j:= \{p\in P\mid
\mu(p)\geq j\}$. Then the homology of $X_{\mathcal{D}}=
\mathbf{hocolim}~\mathcal{D}$ with coefficients in a group
$\mathbb{A}$ admits the decomposition
\begin{equation}
\label{eqn:sumasuma}
 \begin{CD}
 H_\ast(X_\mathcal{D}) &\mbox{ $\cong$ }
& H_\ast(SP^m({X})\times \Delta(P))&
{\displaystyle\bigoplus_{j=m+1}^M}& \!
H_\ast((SP^j({X}), SP^{j-1}({X}))\times \Delta(P_j))\\
&\cong & H_\ast(SP^m({X})\times \Delta(P)) &
{\displaystyle\bigoplus_{p=1}^k}& H_\ast((SP^{c_p}({X}),
SP^{c_{p-1}}({X}))\times \Delta(P_{c_p})).
\end{CD}
\end{equation}
\end{theo}

\medskip\noindent
{\bf Proof:} The result is easily deduced from
Proposition~\ref{prop:injective}. Let
\[
\mathcal{D}_0\preccurlyeq_P\mathcal{D}_1\preccurlyeq_P\ldots\preccurlyeq_P\mathcal{D}_{k-1}
\preccurlyeq_P\mathcal{D}_k = \mathcal{D}
\]
be the sequence of simple diagrams over $P$ where
$\mathcal{D}_j(p) = SP^{c_j}({X})$ if $p\in P_{c_j}$ and
$\mathcal{D}_j(p) = \mathcal{D}(p)$ otherwise, while
$(c_j)_{j=0}^k$ is the sequence defined in (\ref{eqn:assemble}).
In other words $\mathcal{D}_j$ is the simple diagram associated to
the rank function $\mu_j = \mbox{min}\{\mu, c_j\}$. By
Proposition~\ref{prop:injective}
\begin{equation}
\label{eqn:zbir} H_\ast(X_\mathcal{D}) \cong
H_\ast(X_{\mathcal{D}_0})\oplus\bigoplus_{i=1}^k
H_\ast(X_{\mathcal{D}_i}, X_{\mathcal{D}_{i-1}}) .
\end{equation}
Since $\mathcal{D}_0$ is a constant simple diagram we know that
\begin{equation}
\label{eqn:d0} H_\ast(X_{\mathcal{D}_0})\cong
H_\ast(SP^m({X})\times H_\ast(\Delta(P)) .
\end{equation}
Let $\mathcal{E}_i$ and $\mathcal{F}_i$ be constant simple
diagrams over $P_{c_i}$ such that $\mathcal{E}_i(p) =
SP^{c_i}({X})$ and $\mathcal{F}_i(p) = SP^{c_{i-1}}({X})$ for each
$p\in P_{c_i}$. By the excision axiom
 {\small
\begin{equation}
\label{eqn:isecanje} H_\ast(X_{\mathcal{D}_i},
X_{\mathcal{D}_{i-1}})\cong H_\ast(X_{\mathcal{E}_i},
X_{\mathcal{F}_{i}})\cong H_\ast((SP^{c_i}({X}),
SP^{c_{i-1}}({X}))\times \Delta(P_{c_i})) .
\end{equation}
 }
The formulas (\ref{eqn:sumasuma}) follow from (\ref{eqn:zbir})
(\ref{eqn:d0}) (\ref{eqn:isecanje}), the observation that
 {\small
 \begin{equation}
 \label{eqn:potrebno}
 H_\ast((SP^{c_i}({X}), SP^{c_{i-1}}({X})) \times\Delta(P_{c_i}))\cong
\bigoplus_{j=c_{i-1}+1}^{c_i}~H_\ast((SP^{j}({X}),
SP^{j-1}({X}))\times\Delta(P_j))
 \end{equation}}
and the fact that $P_j = P_{c_i}$ for each $j$ in the interval
$(c_{i-1}, c_i]$. \hfill $\square$

\medskip
As an immediate consequence of Proposition~\ref{prop:ample},
Theorem~\ref{thm:suma}, the homotopy equivalence
(\ref{eqn:projection}) and the fact that $\mathcal{E}$ is a simple
diagram, we obtain the following result.

\begin{theo}
\label{thm:vazna} Suppose that $\mathcal{A} = \{F_n^i\}_{i=1}^r$
is a diagram of subspaces of $SP^n({X})$ where $F_n^i = D_i +
SP^{n-\vert D_i\vert}({X})$. Let $P$ be the associated
intersection poset and $\mu : P\rightarrow \mathbb{N}$ the
corresponding rank function. Define $m$ and $M$ respectively as
the minimum and the maximum of the set $\mu(P)\subset \mathbb{N}$
and let $P_j:= \mu^{-1}(\mathbb{N}_{\geq j})$. Then, for the
homology with coefficients in a group $\mathbb{A}$,

\begin{equation}
\label{eqn:vazna}
 \begin{array}{ccl} H_\ast(\bigcup\mathcal{A}) & \cong &
 \, \,  H_\ast(SP^m({X})\times \Delta(P))\, \bigoplus\\
 & &
 {\displaystyle\bigoplus_{j=m+1}^M}~
H_\ast((SP^j({X}), SP^{j-1}({X}))\times \Delta(P_j)).
\end{array}
\end{equation}
\end{theo}

The following corollary of Theorem~\ref{thm:vazna} is needed in
the proof of Theorem~\ref{thm:hend}.

\begin{cor}
\label{cor:hend} Let $\mathcal{A}=\{F_n^i\}_{i=1}^k$ be an
arrangement of subspaces in $SP^n({X})$ where $F_n^i := x_i +
SP^{n-1}({X})$ and $x_1,\ldots, x_k$ are distinct points in ${X}$.
Then, for the homology with rational coefficients,
\begin{equation}
\label{eqn:vazno} \begin{array}{ccl}
 H_\ast(\bigcup\mathcal{A};\mathbb{Q}) &
\cong & \, H_\ast(SP^{n-m}({X}))\otimes
 H_\ast(\Sigma^{m-1}) \, \bigoplus\\ &&
{\displaystyle\bigoplus_{p=0}^{m-2}}~ H_\ast(SP^{n-p-1}({X}),
SP^{n-p-2}({X}))\otimes H_\ast(\Sigma^p)
 \end{array}
\end{equation}
where $\Sigma^p$ is the $p$-skeleton of a simplex $\Sigma$ with
$k$ vertices and $m:= \mbox{\rm min}\{n,k\}$.
\end{cor}

\subsection{Category of simple diagrams}
\label{sec:category}

In this section we take a closer look at the category of simple
diagrams and recast the main results of Section~\ref{sec:homol} in
a form suitable for applications in Section~\ref{sec:AB}. The
emphasis is on functorial properties (naturality) of
decompositions (\ref{eqn:sumasuma}) and (\ref{eqn:vazna}).

\medskip
Let ${\mathcal Pos}$ be the category of finite posets and monotone
(increasing) maps. Let ${\mathcal Rank}$ be the category of
abstract rank functions defined on finite posets. The objects of
$Rank$ are monotone (decreasing) functions $\mu : P\rightarrow
\mathbb{N}$. A monotone map $F : P\rightarrow Q$ defines a
morphism $F : \mu \rightarrow \nu$ of two abstract rank functions
$\mu$ and $\nu$ if $\mu \leq \nu\circ F$, i.e.\ if $\mu(p) \leq
\nu\circ F(p)$ for each $p\in P$. If $F = 1_P$ is the identity map
and $\mu\leq\nu$ then, as in Definition~\ref{def:dominate}, we say
that $\mu$ is dominated by $\nu$ and write $\mu\preccurlyeq\nu$.

\medskip
The category $S$-$Diag$ of simple diagrams is formally isomorphic
to the category $Rank$ of abstract rank functions. Objects of
$S$-$Diag$ are diagrams over finite posets which are simple in the
sense of Definition~\ref{def:simple}. Suppose that $\mathcal{D} :
P\rightarrow Top$ and $\mathcal{E} : Q\rightarrow Top$ are simple
diagrams with the associated rank functions $\mu : P\rightarrow
\mathbb{N}$ and $\nu : Q\rightarrow \mathbb{N}$. Then a morphism
$(F,\alpha) : \mathcal{D}\rightarrow\mathcal{E}$ is defined if $F
: P\rightarrow Q$ is a monotone map and $\mu\leq \nu\circ F$, in
which case $\alpha(p) :
\mathcal{D}(p)\rightarrow\mathcal{E}(F(p))$ is a standard
inclusion in the sense of Definition~\ref{def:standard}. If $F$ is
clear from the context, for example in the case of an identity
map, the corresponding morphism is simply denoted by $\alpha$. A
morphism $(F,\alpha)$ induces a continuous map $X_{(F,\alpha)}  :
X_{\mathcal{D}}\rightarrow X_{\mathcal{E}}$ of the corresponding
homotopy colimits. Again we simplify and often write $\hat\alpha$
instead of $X_{(F,\alpha)}$ if $F : P\rightarrow Q$ is
self-understood.

\medskip
Suppose that $\mu : P\rightarrow \mathbb{N}$ is an object in
$Rank$ and let $\mathcal{D}$ be the associated simple diagram.
Given $j\in \mathbb{N}$, let $\mu_j : P\rightarrow \mathbb{N}$ be
defined by $\mu_j(p) = \mbox{min}\{\mu(p),j\}$ for each $p\in P$.
Define $\mathcal{D}_j$ as the simple diagram associated to $\mu_j$
and let $\alpha_j : \mathcal{D}_j\rightarrow\mathcal{D}$ be the
associated morphism. The associated homotopy colimits
$X_{\mathcal{D}_j}$ are subspaces of $X_{\mathcal{D}}$ which
define a filtration
\begin{equation}
\label{eqn:filtr}
 X_{\mathcal{D}_0} \subseteq X_{\mathcal{D}_1}\subseteq\ldots
 \subseteq X_{\mathcal{D}_j} \subseteq\ldots\subseteq
 X_{\mathcal{D}}.
\end{equation}
Then in light of the decomposition (\ref{eqn:potrebno}),
Theorem~\ref{thm:suma} can be rewritten as follows
\begin{theo}
\label{thm:opetsuma} Assume that $\mathcal{D} : P\rightarrow Top$
is a simple diagram where $\mu : P \rightarrow \mathbb{N}$ is the
corresponding rank function. Let  $P_j:= \{p\in P\mid \mu(p)\geq
j\}$. Then, for the homology with coefficients in an arbitrary
group $\mathbb{A}$,
\begin{equation}
\label{eqn:opetsumasuma}
 H_\ast(X_\mathcal{D};\mathbb{A}) \cong
\bigoplus_{j=0}^\infty H_\ast((SP^j({X}), SP^{j-1}({X}))\times
\Delta(P_j))
\end{equation}
where by definition $SP^{-1}(X) :=  \emptyset =:
\Delta(\emptyset)$.
\end{theo}

The following proposition essentially claims that the
decomposition (\ref{eqn:opetsumasuma}) is natural with respect to
morphisms $(F,\alpha): \mathcal{D}\rightarrow\mathcal{E}$ in the
category $S$-$Diag$.

\begin{prop}
\label{prop:naturalsuma} Suppose that $(F,\alpha):
\mathcal{D}\rightarrow\mathcal{E}$ is a morphism of two simple
diagrams and let $\hat\alpha : X_{\mathcal{D}}\rightarrow
X_{\mathcal{E}}$ be the induced map of the associated homotopy
colimits. If $\{X_{\mathcal{D}_j}\}_{j=0}^{\infty}$ and
$\{X_{\mathcal{E}_j}\}_{j=0}^{\infty}$ are the filtrations  of
$X_{\mathcal{D}}$ and $X_{\mathcal{E}}$  described by {\rm
(\ref{eqn:filtr})} then $\hat\alpha(X_{\mathcal{D}_j})\subseteq
X_{\mathcal{E}_j}$. It follows that there exists a homomorphism
\begin{equation}
\label{eqn:terms}
 H_\ast((SP^j({X}), SP^{j-1}({X}))\times
\Delta(P_j))\longrightarrow H_\ast((SP^j({X}),
SP^{j-1}({X}))\times \Delta(Q_j))
\end{equation}
of the corresponding terms in the decompositions {\rm
(\ref{eqn:opetsumasuma})} of $H_\ast(X_\mathcal{D})$ and
$H_\ast(X_\mathcal{E})$ respectively. The homomorphism {\rm
(\ref{eqn:terms})} is induced by the map $\Delta(F_j):
\Delta(P_j)\rightarrow \Delta(Q_j)$ where $F_j : P_j\rightarrow
Q_j$ is the restriction of $F$ to $P_j$.
\end{prop}

\medskip\noindent
{\bf Proof:} The condition $\mu\leq \nu\circ F$ implies
\[
\mu_j = \mbox{min}\{\mu,j\} \leq \mbox{min}\{\nu\circ F, j\} =
\mbox{min}\{\nu, j\}\circ F = \nu_j\circ F .
\]
It follows that there is a morphism in $S$-$Diag$ of diagrams
$\mathcal{D}_j$ and $\mathcal{E}_j$ and an associated continuous
map $\hat\alpha_j : X_{\mathcal{D}_j} \rightarrow
X_{\mathcal{E}_j}$ of the corresponding homotopy colimits.
Moreover, there is a ladder of commutative diagrams
\begin{equation}
\label{eqn:ladder}
\begin{CD}
\ldots & \longrightarrow & X_{\mathcal{D}_{j-1}} & \longrightarrow
& X_{\mathcal{D}_{j}} & \longrightarrow & X_{\mathcal{D}_{j+1}} &
\longrightarrow & \ldots & \longrightarrow & X_{\mathcal{D}}\\
@. @V{\hat\alpha_{j-1}}VV @V{\hat\alpha_{j}}VV @V{\hat\alpha_{j+1}}VV @. @VV\hat\alpha V\\
\ldots & \longrightarrow & X_{{\mathcal{E}_{j-1}}} &
\longrightarrow & X_{\mathcal{E}_{j}} & \longrightarrow &
X_{\mathcal{E}_{j+1}} & \longrightarrow & \ldots & \longrightarrow
& X_{\mathcal{E}}
\end{CD}
\end{equation}
where the horizontal maps are the inclusions coming from the
filtration (\ref{eqn:filtr}). It follows that there is a
homomorphism
\begin{equation}
\label{eqn:zbirzbir}
 H_\ast(X_{\mathcal{D}}; \mathbb{A}) \cong
 \bigoplus_{j=0}^{\infty} H_{\ast}(X_{\mathcal{D}_j},
 X_{\mathcal{D}_{j-1}}; \mathbb{A}) \longrightarrow
 \bigoplus_{j=0}^{\infty} H_{\ast}(X_{\mathcal{E}_j},
 X_{\mathcal{E}_{j-1}}; \mathbb{A}) \cong
 H_\ast(X_{\mathcal{E}}; \mathbb{A})
\end{equation}
where $X_{\mathcal{D}_{-1}} = X_{\mathcal{D}_{-1}} = \emptyset$.
The final part of Proposition~\ref{prop:naturalsuma} follows from
the naturality of the excision operation or more precisely from
the naturality of the diagram (\ref{eqn:excision}).
\hfill$\square$

\subsection{Homology of the complement $SP^n(M_g)\setminus\bigcup\mathcal{A}$}
 \label{sec:AB}

In this section we focus our attention on the homology of the
complement $SP^n(M_g)\setminus F_n$ of the arrangement
$\mathcal{A}$ where $M_g$ is an orientable surface of genus $g$.
We will be primarily interested in the homology with rational
coefficients. By Poincar\' e duality, the evaluation of these
groups is equivalent to the evaluation of the homology of the pair
$H_\ast((SP^n(M_g), F_n);\mathbb{Q})$ which is directly related to
the evaluation of the kernel $A_j$ and the cokernel $B_j$ of the
homomorphism
\begin{equation}
\label{eqn:kercoker} H_j(\bigcup\mathcal{A}) \longrightarrow
H_j(SP^n(X))
\end{equation}
for an arrangement $\mathcal{A}$ in $SP^n(X)$ where $X$ is an {\em
admissible} space.

\medskip
Given a map of diagrams $\alpha :
\mathcal{E}_1\rightarrow\mathcal{E}_2$, there is a commutative
square
\begin{equation}
\label{eqn:square}
\begin{CD}
\mbox{\bf hocolim}~\mathcal{E}_1 @>>> \mbox{\bf
hocolim}~\mathcal{E}_2\\
@VVV @VVV \\
\mbox{\bf colim}~\mathcal{E}_1 @>>> \mbox{\bf colim}~\mathcal{E}_2
\end{CD}
\end{equation}
In light of the fact that for each arrangement of subspaces
$\mathcal{F}$ and the corresponding diagram of spaces
$\mathcal{D}_{\mathcal{F}}$
\[
\bigcup~\mathcal{F} \cong \mbox{\bf
colim}~\mathcal{D}_{\mathcal{F}},
\]
the square (\ref{eqn:square}) allows us to compare spaces $F_n =
\bigcup\mathcal{A}$ and $SP^n({X})$ by comparing the associated
diagrams.

The diagrams $\mathcal{C},\mathcal{D}, \mathcal{E}$, the
associated intersection poset $P$ etc.\ have the same meaning as
in Section~\ref{sec:homol}. Let $\hat{P} := P\cup\{0\}$ be the
poset $P$ with added a possibly new minimum element $0$. Let
$\mathcal{D}_0 : \hat{P}\rightarrow Top$ be a diagram of spaces
and inclusion maps defined by $\mathcal{D}_0(0) = SP^n({X})$ and
$\mathcal{D}_0(p) = \mathcal{D}(p)$ for each $p\in P$. Similarly,
$\mathcal{C}_0$ and $\mathcal{E}_0$ are defined as the diagrams
over $\hat{P}$ which extend the diagrams $\mathcal{C}$ and
$\mathcal{E}$ respectively, such that
$\mathcal{C}_0(0)=\mathcal{E}_3(0) = SP^n({X})$. The inclusion $e
: P\rightarrow \hat{P}$ extends to the corresponding morphisms of
diagrams
\[
{\mathcal{C}}\stackrel{\gamma}{\longrightarrow}{\mathcal{C}_3}
 \qquad
{\mathcal{D}}\stackrel{\delta}{\longrightarrow}{\mathcal{D}_0}
 \qquad
{\mathcal{E}}\stackrel{\eta}{\longrightarrow}{\mathcal{E}_1}
\]
Moreover, there is a commutative diagram
\begin{equation}
\label{eqn:kvadrati}
\begin{CD}
 X_\mathcal{D} @>\hat\alpha>>
X_\mathcal{C} @<\hat\beta << X_\mathcal{E}
\\
@V\hat\delta VV @V\hat\gamma VV @V\hat\eta VV\\
X_{\mathcal{D}_0}@>\hat{\alpha}_0>>
X_{\mathcal{C}_0}@<\hat{\beta}_0<< X_{\mathcal{E}_0}
\end{CD}
\end{equation}
where all horizontal maps are homotopy equivalences and as before
$X_{\mathcal{G}} = \mbox{\bf hocolim}~\mathcal{G}$. An instance of
the commutative diagram (\ref{eqn:square}) is the following square
\begin{equation}
\begin{CD}
X_\mathcal{D} @>\hat\delta>> X_{\mathcal{D}_0} \\
@VVV @VVV \\
\bigcup\mathcal{A} @>>> SP^n(M_g)
\end{CD}
\end{equation} By Projection Lemma, \cite{WZZ99} Lemma~4.5, the
vertical arrows in this diagram are homotopy equivalences. In
light of (\ref{eqn:kvadrati}) we conclude that the groups $A_d$
and $B_d$ are respectively the kernel and the cokernel of the map
\begin{equation}
\label{eqn:extended}
\begin{CD}
H_d(X_\mathcal{E}) @>\eta_\ast >> H_d(X_{\mathcal{E}_0}) .
\end{CD}
\end{equation} By Theorem~\ref{thm:opetsuma} each of the
groups $H_d(X_\mathcal{E})$ and $H_d(X_{\mathcal{E}_0})$ admits a
direct sum decomposition of the form (\ref{eqn:opetsumasuma}). By
Proposition~\ref{prop:naturalsuma} this decomposition is natural
and the corresponding terms are mapped to each other. More
precisely there is a homomorphism
\begin{equation}
H_\ast((SP^j({X}), SP^{j-1}({X}))\times \Delta(P_j))
\longrightarrow H_\ast((SP^j({X}), SP^{j-1}({X}))\times
\Delta(\hat{P}_j))
\end{equation}
induced by the map $\Delta(e_j) :
\Delta(P_j)\rightarrow\Delta(\hat{P}_j)$ where $e_j$ is an
inclusion map.

Note that $\Delta(\hat{P}_j)$ is contractible whenever
$\hat{P}_j\neq\emptyset$. It immediately follows that for the
homology with {\em rational coefficients}, $A_d\cong
A_d^{(1)}\oplus A_d^{(2)}$ where
\begin{equation}
A_d^{(1)}\cong \bigoplus_{j=0}^{+\infty} \bigoplus_{p+q = d}^{q >
0} H_p(SP^j({X}), SP^{j-1}(X))\otimes H_q(\Delta({P}_j))
\end{equation}
\begin{equation}
A_d^{(2)}\cong \bigoplus_{j=9}^{\infty} H_d(SP^j({X}),
SP^{j-1})\otimes \Omega_j
\end{equation}
where $\Omega_j := {\rm Ker}\{H_0(\Delta(P_j))\rightarrow
H_0(\Delta(\hat{P}_j))\} = \tilde{H}_0(\Delta(P_j))$ is the
reduced, $0$-dimensional homology of $P_j$. Similarly,
\begin{equation}
B_d \cong H_d(SP^n(X),SP^M({X}))
\end{equation}
where $M={\rm max}\{\mu(p)\}_{p\in P}$ is the maximum of the
associated rank function.

\subsection{An auxiliary result}
\label{sec:main}

\begin{prop}
\label{prop:main} Suppose that $X$ is an {\em admissible} space in
the sense of Definition~\ref{def:admissible}. Given a finite set
$F =\{y,y_1,\ldots, y_d\}$ in $Z$, let $C$ be an associated
contractible superset of $F$. Let $Y = y_1+y_2+\ldots + y_d \in
SP^d(X)$. Then the inclusion map
\begin{equation}
\label{eqn:inclusion2} Y + SP^{n-d}(X)
\stackrel{\alpha}{\longrightarrow} SP^d(C) + SP^{n-d}(X)
\end{equation}
is a homotopy equivalence.
\end{prop}

\medskip\noindent
{\bf Proof:} We prove the lemma first in the special case when all
points $y_i$ coincide. More precisely we prove that the map $d\, y
+ SP^{n-d}(X) \hookrightarrow SP^d(C) + SP^{n-d}(X)$ is a homotopy
equivalence.

It follows from the assumptions on $X, C$ and $y\in C$ that there
exists a homotopy $h : X\times I\rightarrow X$, keeping $y$ fixed
and $C$ invariant, such that $p : X\rightarrow X$, defined by
$p(x):= h(x,1)$ satisfies the condition $p(C)=\{y\}$ and $1_X =
h(\cdot, 0)$ is the identity map.

Let $H : SP^n(X)\times I \rightarrow SP^n(X)$ be the homotopy
induced on $SP^n(X)$ by $h$. In other words if $Z = z_1+\ldots +
z_n\in SP^n(X)$, then $H(Z,t) := h(z_1,t) +\ldots +h(z_n,t)$. We
observe that both $V:= d\, y + SP^{n-d}(X)$ and $W:= SP^d(C) +
SP^{n-d}(X)$ are $H$-invariant, so the restrictions of $H$ on
these subspaces define the homotopies $H^1 : V\times I\rightarrow
V$ and $H^2: W\times I\rightarrow W$. The map $p :X\rightarrow X$
induces a map from $SP^n(X)$ to $SP^n(X)$ which restrict to a map
$\hat{p} : W\rightarrow V$. It turns out that $\hat{p}$ is a
homotopy inverse to the inclusion $i : V\hookrightarrow W$.
Indeed,
\[
H^1 : 1_V \simeq \hat{p}\circ i \qquad \mbox{ and } \qquad H^2 :
1_W\simeq i\circ \hat{p} \, .
\]

Now we turn to the case of a general divisor $Y = y_1+\ldots
+y_d\in SP^d(X)$ and the associated inclusion map $Y + SP^{n-d}(X)
\stackrel{\alpha}{\longrightarrow} SP^d(C) + SP^{n-d}(X)$. Let
$\phi : d\, y+ SP^{n-d}(X)\rightarrow Y + SP^{n-d}(X)$ be the
homeomorphism defined by $\phi (d\, a + Z) = Y + Z$.

The key observation is the following equality
\[
\hat{p}\circ\alpha\circ\phi = \hat{p}\circ i .
\]

Since $\hat{p}$ and $i$ are homotopy equivalences and $\phi$ is a
homeomorphism, we conclude that $\alpha$ is also a homotopy
equivalence. \hfill $\square$

\medskip

\begin{defin}
\label{def:admissible} A space $X$ is called {\em admissible} if
for each finite collection of points $F = \{y,x_1,\ldots
,x_k\}\subset X$, there exists a subspace $C\subset X$ such that
\begin{itemize}
 \item[(a)]  $C$ contains $F$ as a subset,
 \item[(b)]  $C$ can be continuously deformed to $y$ inside $C$ keeping the point $y$
 fixed, i.e.\ $1_C \simeq c_y \, (\mbox{\rm rel } y)$ where $c_y(x) =
 y  \mbox{ {\rm for each} } x\in C$,
 \item[(c)] the inclusion map $i : C\hookrightarrow X$ is a closed
 cofibration.
\end{itemize}
\end{defin}

\begin{rem}
{\rm All connected spaces that can be triangulated are {\em
admissible} in the sense of Definition~\ref{def:admissible}. This
is a very large class including all connected $CW$-complexes or
connected semi-algebraic sets.}
\end{rem}

\section{Applications}
\label{sec:applications}

\subsection{End spaces}

\begin{defin}
\label{def:appl} Suppose that $Y$ is a locally compact, Hausdorff
space. Let ${\mathcal P} = ({\mathcal K},\subseteq )$ be the poset
of all compact subspaces of $Y$. Let ${\mathcal D} : {\mathcal
P}\rightarrow Top$ be the diagram of topological spaces over
$\mathcal P$ defined by ${\mathcal D}(K) := Y\setminus K$ for
$Q\in\mathcal M$. The {\em end space} $e(Y)$ of $Y$ is by
definition the homotopy limit of the diagram ${\mathcal D}$,
  \[ e(Y):= \mbox{ \bf holim } {\mathcal D} .\]
\end{defin}
The reference \cite{ends} is recommended as a valuable source of
information about the general theory, the history, and some of the
latest applications of end spaces. The end space is obviously a
topological invariant of $Y$, that is if $Y$ and $Y'$ are two
locally compact, Hausdorff spaces such that $Y\cong Y'$ then the
associated end spaces $e(Y)$ and $e(Y')$ are also homeomorphic.
Consequently the homotopy type of $e(Y)$, its homology etc. are
homeomorphism invariants of $Y$. If $Y$ admits a cofinal sequence
\[
K_0\subseteq K_1 \subseteq \ldots \subseteq K_m \subseteq \ldots
\]
of compact sets in $Y$, i.e.\ a sequence such that $Y =
\bigcup_{m=0}^{\infty}~K_m$, then
\[ e(Y) \simeq \mbox{\bf holim}_{m\mapsto\infty} Y\setminus K_m .\]

The following proposition allows us to ``compute'' the end space
in the case the inclusion map $K_m\hookrightarrow K_{m+1}$ is a
homotopy equivalence for each $m$.

\begin{prop} Let ${\mathcal E} : \mathbb{N}\rightarrow Top$ be a
diagram of topological spaces over $\mathbb{N}$ such that
${\mathcal D}(m) \rightarrow {\mathcal D}(m+1)$ is a homotopy
equivalence for each $m$. Then
\[
\mbox{ \bf holim } {\mathcal E} \simeq {\mathcal D}(0).
\]
\end{prop}

\medskip
As an illustration here is a computation of the (stable) homotopy
type of the end space of some interesting spaces.

\begin{prop}
Let $Y = S^n\setminus X$ where $X$ is a closed set in the sphere
$S^n$. Let $D_n(A)$ be the stable homotopy type of the geometric
dual of $A\subset S^n$, {\rm \cite{SpaWhi55}}. Then the end space
of $Y$ is stably equivalent to the geometric dual of the disjoint
sum $X\cup D_n(X)$,
\[
e(Y) \simeq_P D_n(X\cup D_n(X)) \simeq_P S^{n-1}\vee X \vee
D_n(X).
\]
\end{prop}
\noindent{\bf Proof:} In light of the formula
 \[
 D_n(A\cup B) \simeq_P D_n(S^0\vee A\vee B) \simeq_P S^{n-1}\vee
 D_n(A)\vee D_n(B)
 \]
where $A$ and $B$ are disjoint, compact subsets of $S^n$, it is
sufficient to observe that compact sets in $Y = S^n\setminus X$
which are deformation retracts of $Y$ are cofinal in the poset of
all compact subsets of $Y$.

\subsection{End (co)homology groups} \label{sec:hend}

The diagram ${\mathcal D} : {\mathcal P}\rightarrow Top$ from
Definition~\ref{def:appl},  in combination with standard functors,
provides other interesting invariants or the homeomorphism type of
$Y$.

\begin{defin}
\label{def-homology} Let $\phi : Top \rightarrow Ab$ be a
covariant (contravariant) functor from spaces to Abelian groups.
Let $\phi({\mathcal D}) : {\mathcal P} \rightarrow Ab$ be the
associated (co)diagram defined by the correspondence $K\mapsto
\phi(Y\setminus K)$. The (co)limit of $\phi({\mathcal D})$ is an
Abelian group which is called the end $\phi$-group associated to
$Y$. If $\phi$ is a homology (cohomology) functor with rational
coefficients, then $\phi({\mathcal D})$ is denoted by ${\mathcal
H}_\ast$ and ${\mathcal H}^\ast$ respectively. In other words
${\mathcal H}_{\ast}(K) := H_{\ast}(Y\setminus K; \mathbb{Q})$ and
${\mathcal H}^{\ast}(K) := H^{\ast}(Y\setminus K; \mathbb{Q})$.
The associated {\em end (co)homology groups} are defined by

\[
E_{\ast}(Y) := \mbox{ \bf lim }{\mathcal H}_\ast , \qquad
E^{\ast}(Y) := \mbox{ \bf colim }{\mathcal H}^\ast .
\]
\end{defin}

In the following theorem we compute the group
$E^{\ast}(SP^n(M_{g,k}))$ where $M_{g,k}$ is the orientable
surface of genus $g$ punctured at $k$ distinct points. Sometimes
we call the surface $M_{g,k}$ a $(g,k)$-amoeba, see
Figure~\ref{fig:amoeba} where the amoebas $M_{1,3}$ and $M_{2,1}$
are shown. Note that although the amoebas $M_{g,k}$ and
$M_{g',k'}$ are homeomorphic if and only if $(g,k)=(g',k')$, they
have the same homotopy type if and only if $2g+k=2g'+k'$.

\begin{theo}
\label{thm:hend} Let $M_{g,k} = M_g\setminus\{x_1,\ldots, x_k\}$
be a $(g,k)$-amoeba, i.e.\ the Riemann surface of genus $g$ with
$k$ distinct points removed. Let $SP^n(M_{g,k})$ be the associated
symmetric product. If
\[
E^{p}(SP^n(M_{g,k})) := \mbox{ \bf colim }{\mathcal H}^{p} =
\mbox{ \bf colim}_{W\in {\mathcal P}} H^{p}(QP^n(M_{g,k})\setminus
K; \mathbb{Q})
\]
is the associated $p$-dimensional end cohomology group then
\begin{equation}
\label{eqn:glavnica}  {\rm rank}(E^{p}(SP^n(M_{g,k}))) =
\left\{
\begin{array}{ll}
 \binom{2g+k-1}{p}, & p\leq n-2\\
 \binom{2g+k}{n} - \binom{2g}{n}, & p=n-1 \mbox{ {\rm or} } p=n\\
 \binom{2g+k-1}{2n-1-p}, & p\geq n+2 \, .
\end{array} \right.
\end{equation}
\end{theo}

\medskip
\noindent{\bf Proof:} Let us choose a local metric in the vicinity
of each point $x_i$, say by choosing a local coordinate system in
the neighborhood of each of these points. Let $V^i_m =
V(x_i,\frac{1}{m})$ be an open disc in $M_g$ with the center at
$x_i$ of radius $\frac{1}{m}$, defined relative to the chosen
metric. Let ${W}^i_m:= V^i_m\setminus\{x_i\}$ be the corresponding
punctured disc. Then $C_m = M_g\setminus
\bigcup_{i=1}^k~V(x_i,\frac{1}{m})$ is a compact subset in
$M_{g,k}$ and $K_m := SP^n(C_m)$ is a compact subset in
$SP^n(M_{g,k})$. The sequence $\{K_m\}_{m=1}^\infty$ is cofinal in
the poset ${\mathcal P}$ of all compact subsets in $Y =
SP^n(M_{g,k})$ since $\bigcup_{m=1}^{+\infty}~K_m =
SP^n(M_{g,k})$. Note that $Y\setminus K_m$ is described as the
space of all divisors $D\in SP^n(M_{g,k})$ such that $D\cap
V^i_m\neq\emptyset$ for some $i$.

\medskip\noindent
{\bf Claim~1:} The inclusion $Y\setminus K_{m+1}\hookrightarrow
Y\setminus K_m$ induces an isomorphism $H^{\ast}(Y\setminus K_m;
\mathbb{Q}) \rightarrow H^{\ast}(Y\setminus K_{m+1}; \mathbb{Q})$
of the associated cohomology groups.

\medskip\noindent
The claim is an easy consequence of Poincar\' e duality. Indeed,
both $Y$ and its compact subsets $K_m$ and $K_{m+1}$ are
identified with the corresponding subsets of the manifold
$SP^n(M_g)$. Let $F_n:= SP^n(M_g)\setminus Y$ be the subspace of
all divisors $D\in SP^n(M_g)$ such that $D\cap\{x_1,\ldots,
x_k\}\neq\emptyset$. Equivalently,
\[
F_n := \bigcup_{j=1}^k (x_j + SP^{n-1}(M_g)) .
\]
There is a duality isomorphism
\begin{equation}
\label{eqn:vazna2} H^{\ast}(Y\setminus K_j) \rightarrow
H_{2n-\ast}(SP^n(M_g); F_n\cup K_j)
\end{equation}
natural with respect the inclusions $K_j\hookrightarrow K_{j+1}$.
Then the claim follows from the five lemma and the fact that $K_m$
is a deformation retract of $K_{m+1}$, which in turn implies the
isomorphism $H_{\ast}(F_n\cup K_m)\rightarrow H_{\ast}(F_n\cup
K_{m+1})$.

\medskip\noindent
A consequence of the claim is the isomorphism
\begin{equation}
\label{eqn:iso}
 E^{p}(SP^n(M_{g,k})) \cong H^{p}(Y\setminus K_m) \cong
H_{2n-p}(SP^n(M_g); F_n\cup K_m).
\end{equation}
Hence, in order to complete the proof of the theorem it is
sufficient to compute the group  $H_{d}(SP^n(M_g); F_n\cup K_m)$.

\begin{prop}
\label{prop:linklink} Let $\Delta_d^{n,g}$ be the rank of the
group $H_d(SP^n(M_g); F_n\cup K_m)\cong
H^{2n-d}(SP^n(M_g)\setminus F_n\cup K_m)$. Then,

\begin{equation}
\label{eqn:osnova} \Delta_d^{n,g} = \left\{ \begin{array}{ll}
 \binom{2g+k-1}{d-1}, & d\leq n-1\\
 \binom{2g+k}{n} - \binom{2g}{n}, & d=n \mbox{ {\rm or} } d=n+1\\
 \binom{2g+k-1}{2n-d}, & d\geq n+2 \, .
\end{array} \right.
\end{equation}
\end{prop}

\medskip\noindent
{\bf Proof:} The spaces $F_n$ and $K_m$ are disjoint compact
subspaces of $SP^n(M_g)$ and $K_m \simeq SP^n(M_{g,k})$ has the
homotopy type of an $n$-dimensional CW-complex.  Consequently, the
exact sequence of the pair $(SP^n(M_g), F_n\cup K_m)$ has the form
 {\small
\begin{equation}
\label{eq:MV1}
\begin{array}{cccccl}
  H_d(F_n)\oplus H_d(K_m) &
 \stackrel{\Lambda_d}{\longrightarrow}&  H_{d}(SP^n(M_g)) & \longrightarrow
  & H_{d}(SP^n(M_g); F_n\cup
K_m)  & \longrightarrow \\
 H_{d-1}(F_n)\oplus H_{d-1}(K_m) &
 \stackrel{\Lambda_{d-1}}{\longrightarrow} &  H_{d-1}(SP^n(M_g)) & \longrightarrow  & \ldots
 &
\end{array}
\end{equation}
 }
An immediate consequence of (\ref{eq:MV1}) is the isomorphism
\[
 H_{d}(SP^n(M_g); F_n\cup K_m) \cong A_{d-1}\oplus B_d
\]
where $A_d := {\rm Ker}(\Lambda_d)$ and $B_d := \mbox{\rm
Coker}(\Lambda_d)$. Note that $\Lambda_d = \alpha_d + \beta_d$
where $\alpha_d : H_d(F_n)\rightarrow H_d(SP^n(M_g))$ and $\beta_d
: H_d(SP^n(M_{g,k}))\rightarrow H_d(SP^n(M_g))$ are the
homomorphisms induced by the inclusions $F_n\hookrightarrow
SP^n(M_g)$ and $SP^n(M_{g,k})\hookrightarrow SP^n(M_g)$.

\medskip\noindent
{\bf Claim~2}: Suppose that $\alpha : A\rightarrow C$ and $\beta :
B\rightarrow C$ a linear maps of vector spaces.  If $\alpha +\beta
: A\oplus B\rightarrow C$ is the map defined by $(\alpha +
\beta)(x\oplus y) := \alpha(x) + \beta(y)$ then
\begin{equation}
\label{eqn:claim2}
 \vert {\rm Ker}(\alpha +\beta)\vert = \vert {\rm Ker}(\alpha)\vert +
 \vert {\rm Ker}(\beta)\vert + \vert{\rm Im}(\alpha)\cap {\rm
 Im}(\beta)\vert
\end{equation}
where $\vert V\vert$ is the dimension of $V$.

\medskip
An immediate consequence of the Claim~2 is the following equation
\begin{equation}
\label{eqn:prima}
 \begin{array}{ccl}
\Delta_d^{n,g} & = &\vert{\rm Ker}(\alpha_{d-1})\vert + \vert{\rm
Ker}(\beta_{d-1})\vert + \vert{\rm
Im}(\alpha_{d-1}\cap{\rm Im}(\beta_{d-1})\vert + \\
& + & \Theta_d^{n,g} -\vert{\rm Im}(\alpha_d)\vert - \vert{\rm
Im}(\beta_d)\vert + \vert{\rm Im}(\alpha_d)\cap {\rm
Im}(\beta_d)\vert
\end{array}
\end{equation}
where $\Delta_d^{n,g}:= \vert A_{d-1}\vert + \vert B_d\vert$ and
$\Theta_d^{n,g} := \vert H_d(SP^n(M_g))\vert$. The long exact
sequence of the pair $(SP^n(M_g), F_n)$ and the Poincar\' e
duality imply that
\begin{equation}
\label{eqn:sekunda}
 \begin{array}{cclc} \vert {\rm Ker}(\alpha_{d-1})\vert +
\Theta_d^{n,g} - \vert{\rm Im}(\alpha_d)\vert & = &  \vert
H_d(SP^n(M_g); F_n)\vert & = \\
 \vert H^{2n-d}(SP^n(M_g)\setminus
F_n)\vert & = & \Phi_d^{n,g,k} &
\end{array}
\end{equation}
where $\Phi_d^{n,g,k} = \binom{2g+k-1}{2n-d}$ if $n\leq d$ and $0$
otherwise. It follows from equation (\ref{eqn:sekunda}),
Proposition~\ref{prop:binom1} and Proposition~\ref{prop:binom2}
that in order to compute the quantity $\Delta_d^{n,g}$ we need to
discuss the four cases $d\leq n-1, d=n, d=n+1$ and $d\geq n+2$.
The rest is an elementary calculation. This completes the proof of
both Proposition~\ref{prop:linklink} and Theorem~\ref{thm:hend}.
\hfill $\square$

\subsection{Commutative $(m+k,m)$-groups}
\label{sec:commutative}

A  commutative $(m+k,m)$-groupoid is a pair $(X,\mu)$ where the
``multiplication'' $\mu$ is a map $\mu : SP^{m+k}(X)\rightarrow
SP^m(X)$. The operation $\mu$ is associative if for each $c\in
SP^{m+2k}(X)$ and each representation $c = a + b$, where $a\in
SP^{m+k}(X)$ and $b\in SP^k(X)$, the result $\mu(\mu(a)\ast b)$ is
always the same, i.e.\ independent from the particular choice of
$a$ and $b$ in the representation $c = a + b$. A commutative and
associative $(m+k,m)$-groupoid is a $(m+k,m)$-group if the
equation $\mu(x + a) = b$ has a solution $x\in SP^m(X)$ for each
$a\in SP^k(X)$ and $b\in SP^m(X)$. The $(2,1)$-groups are
essentially the groups in the usual sense of the word. If $X$ is a
topological space then $(X,\mu )$ is a topological $(m+k,m)$-group
if it is a $(m+k,m)$-group and the map $\mu :
SP^{m+k}(X)\rightarrow SP^m(X)$ is continuous.

For the motivation and other information about commutative
$(m+k,m)$-groups the reader is referred to \cite{TD2001},
\cite{TD92}. The only known surfaces that support the structure of
a $(m+k,m)$-group for $(m+k,m)\neq (2,1)$ are of the form
$\mathbb{C}\setminus A$ where $A$ is a finite set. It was proved
in \cite{TD92}, see also Theorem~6.1 in \cite{TD2001} that if
$(M,\mu)$ is a locally Euclidean, topological, commutative,
$(m+k,m)$-group then $M$ must be an orientable $2$-manifold.
Moreover, a $2$-manifold that admits the structure of a
commutative $(m+k,m)$-group satisfies a strong necessary condition
that the symmetric power $SP^m(M):= M^m/S_m$ is of the form
$\mathbb{R}^u\times (S^1)^v$.

This was a motivation for the authors to formulate in
\cite{BGZ-Israel} the following problems

\begin{enumerate}
\item[(A)]
 To what extent is the topology of a surface $M$
determined by the topology of its symmetric product $SP^m(M)$ for
a given $m$?

\item[(B)] Are there examples of non-homeomorphic (open) surfaces
$M$ and $N$ such that the associated symmetric products $SP^m(M)$
and $SP^m(N)$ are homeomorphic?
\end{enumerate}

\begin{figure}
\centering
\includegraphics[scale=0.80]{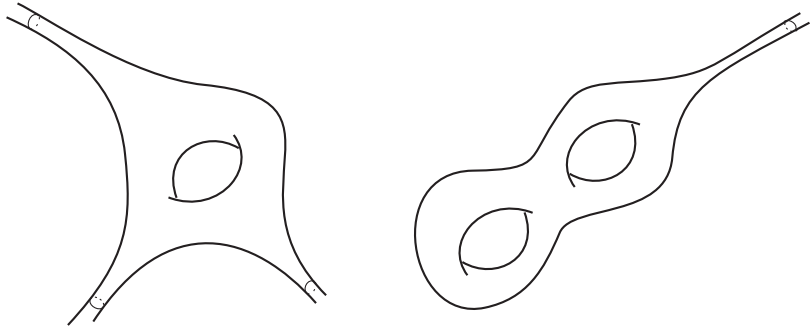}
\caption{$M_{1,3}$ and $M_{2,1}$} \label{fig:amoeba}
\end{figure}

These problems are particularly interesting for the so called
$(g,k)$-amoebas $M_{g,k}$, the surfaces defined by  $M_{g,k} :=
M_g\setminus\{x_1,\ldots, x_k\}$ where $M_g$ is the Riemann
surface of genus $g$, see Figure~\ref{fig:amoeba}.

 Since the end homology groups are
invariants of homeomorphism types, a consequence of
Theorem~\ref{thm:hend} is that both the genus $g$ and the number
$k$ of points removed can be recovered from the knowledge of the
end space $e(SP^n(M_{g,k}))$ in the case $2g\geq n$ (question
(A)). As a corollary we obtain the main result of
\cite{BGZ-Israel} (Theorem~1.1) which says that there exist open,
orientable surfaces $M$ and $N$ such that the associated symmetric
products $SP^{m}(M)$ and $SP^{m}(N)$ are not homeomorphic although
they have the same homotopy type. More precisely, this is always
true if $M = M_{g,k}$ and $N = M_{g',k'}$ {\rm ($k,k'\geq 1$)} and
\begin{itemize}
\item $2g + k = 2g' + k'$, \item $g\neq g'$ and ${\rm
max}\{g,g'\}\geq m/2$ .
\end{itemize}
This result puts some restrictions on potential examples asked for
in question (B). Note that this result was proved in
\cite{BGZ-Israel} by completely different methods based on the
evaluation of the signature of general symmetric products
$SP_G(M_{g,k}):= (M_{g,k})^n/G$ where $G$ is a subgroup of $S_n$.
The fact that two different methods led to the exactly same
necessary condition $2g\geq n$ is intriguing and may be an
indication that the answer to the question (B) is in fact
positive. This however remains an interesting open problem which
can be restated as follows.

\medskip
\noindent {\bf Question:} Could it be that some symmetric powers
of amoebas $M_{1,3}$ and $M_{2,1}$ or more generally amoebas
$M_{g,k}$ and $M_{g',k'}$ where $(g,k)\neq (g',k')$ but
$2g+k=2g'+k'$ are actually homeomorphic?

\section{Appendix}

\subsection{Homology of symmetric products}
\label{sec:homology}

It is well known \cite{Milgram} \cite{Kallel_1998} that the
Pontriagin algebra of $SP^{\infty}(M_g)$ has the form

\begin{equation}
\label{eqn:Pontriagin1} H_{\ast}(SP^{\infty}(M_g);\mathbb{A})
\cong \Lambda(e_1,e_2,\ldots, e_{2g})\otimes \Gamma[M]
\end{equation}
where $\Lambda(e_1,\ldots, e_{2g})$ is an exterior algebra and
$\Gamma[M]$ is the divided power algebra with generators
$\gamma_k, \, k=1,2,\ldots \, .$ On tensoring with $\mathbb{Q}$
one has an isomorphism
\begin{equation}
\label{Pontriagin2} H_{\ast}(SP^{\infty}(M_g);\mathbb{Q}) \cong
\Lambda(e_1,e_2,\ldots, e_{2g})\otimes \mathbb{Q}[\gamma]
\end{equation}
where $\mathbb{Q}[\gamma]$ is a polynomial algebra. The following
classes
\[
E_{(I,q)} =  e_I \gamma^q = e_{i_1}e_{i_2}\ldots
e_{i_p}\cdot\gamma^q, \, I = \{i_1<i_2<\ldots <i_p\}
\]
form an additive basis of (\ref{Pontriagin2}). Moreover,
\begin{equation}
\label{eqn:Steen1} H_{\ast}(SP^{\infty}(M_g)) \cong
\bigoplus_{m\in \mathbb{N}} H_{\ast}(SP^m(M_g), SP^{m-1}(M_g))
\end{equation}
and $H_{\ast}(SP^m(M_g))$ is the subgroup spanned by all classes
$E_{(I,q)} =  e_{i_1}e_{i_2}\ldots e_{i_p}\gamma^q$ where $p+q\leq
m$.

\medskip
Recall \cite{Milgram} \cite{Kallel_1998} that by a classical
result of Steenrod, the decomposition (\ref{eqn:Steen1}) holds for
all connected $CW$-complexes $X$, in particular $H_\ast(SP^m(X);
{\mathbb{A}}) \rightarrow H_\ast(SP^n(X); {\mathbb{A}})$ is always
a monomorphism for $m\leq n$. Similarly, the isomorphism
(\ref{Pontriagin2}) can be seen as an instance of the celebrated
result of Dold and Thom \cite{DoldThom} which says that infinite
symmetric products admit a decomposition into a product of
Eilenberg-Mac~Lane spaces
\begin{equation}
\label{eqn:DoThom} SP^\infty(X) \simeq \prod_{\nu\geq 0}~
K(\tilde{H}_\nu; \mathbb{A}), \nu) .
\end{equation}

\medskip
Suppose that $Y = \bigvee_{j=1}^{m}~S^1_{j}$ is a wedge of $m$
circles. Then it is not difficult to show directly that
\begin{equation}
H_\ast(SP^{\infty}(Y;\mathbb{Q})) \cong \Lambda(e_1,\ldots, e_m) .
\end{equation}
Since $M_{g,k} = M_g\setminus\{x_1,\ldots, x_k\} \simeq
\bigvee_{j=1}^{2g+k-1}~S^1_j$ we conclude that
\begin{equation}
\label{Pontriagin3} H_{\ast}(SP^{\infty}(M_{g,k});\mathbb{A})
\cong \Lambda(e_1,e_2,\ldots, e_{2g+k-1}) .
\end{equation}
The group $H_{\ast}(SP^m(M_{g,k});\mathbb{Q})$ is generated by the
classes $F_I = e_{i_1}\ldots e_{i_p}$ where $I = \{i_1,\ldots,
i_p\}\subset [2g+k-1]$ and $p\leq m$. Moreover, the generators can
be chosen so that $\beta_\ast(F_I) =E_{I,0} = e_I$ if $I\subset
[2g]$ and $0$ otherwise where $\beta_\ast :
H_\ast(SP^m(M_{g,k});\mathbb{Q})\rightarrow
H_\ast(SP^n(M_g);\mathbb{Q})$ is the map induced by the inclusion
$SP^m(M_{g,k})\hookrightarrow SP^m(M_g)$. As an immediate
consequence we have the following result needed in the proof of
Theorem~\ref{thm:hend} in Section~\ref{sec:applications}. As
before, $\vert V\vert$ is the rank of a vector space $V$.

\begin{prop}\label{prop:binom1}
\begin{equation}
\begin{array}{ccl}
 \vert{\rm Ker}(\beta_d)\vert & = & \left\{
\begin{array}{cl} \binom{2g+k-1}{d} - \binom{2g}{d}, & d\leq n \\
              0, & d \geq n+1
\end{array} \right. \\
\vert{\rm Im}(\beta_d)\vert & = & \left\{
\begin{array}{cl} \binom{2g}{d} , & d\leq n \\
              0, & d \geq n+1
\end{array} \right.
\end{array}
\end{equation}
\end{prop}

As a consequence of results from Section~\ref{sec:category} we
deduce the following proposition which is also used in the proof
of Theorem~\ref{thm:hend}.

\begin{prop}
\label{prop:binom2} Let $\alpha_p : H_p(F_n)\rightarrow
H_p(SP^n(M_g))$ and $\beta_p : H_p(SP^n(M_{g,k}))\rightarrow
H_p(SP^n(M_g))$ be the homomorphisms induced by the inclusions
$F_n\hookrightarrow SP^n$ and $SP^n(M_{g,k})\hookrightarrow
SP^n(M_g)$. If $\alpha_p + \beta_p : H_p(F_n)\oplus
H_p(SP^n(M_{g,k})\rightarrow H_p(SP^n(M_g))$ is the map defined by
$(\alpha_p+\beta_p)(x\oplus y) = \alpha_p(x) + \beta_p(y)$, then
\begin{equation}
\vert {\rm Im}(\alpha_p) \cap {\rm Im}(\beta_p)\vert  =  \left\{
\begin{array}{cl} \binom{2g}{p} , & p\leq n-1 \\
              0, & p \geq n
\end{array} \right.
\end{equation}
\end{prop}

\medskip\noindent
{\bf Proof:} By the results of Sections~\ref{sec:category} and
\ref{sec:AB} the group $H_p(F_n)$ has the following decomposition
\begin{equation}
\label{eqn:decomp} H_p(F_n) \cong \bigoplus_{\nu=0}^{+\infty}
\bigoplus_{i+j=p} H_i(SP^\nu(M_g), SP^{\nu -1}(M_g))\otimes
H_j(\Delta(P_\nu)).
\end{equation}
Define
\[
H_p(F_n)^{(0)} := \bigoplus_{\nu=0}^{+\infty}~ H_p(SP^\nu(M_g),
SP^{\nu -1}(M_g))\otimes H_0(\Delta(P_\nu))
\]
and let $\alpha_p' : H_p(F_n)^{(0)}\rightarrow H_p(SP^n(M_g))$ be
the restriction of $\alpha_p$ on $H_p(F_n)^{(0)}$. Let us observe
that
\begin{equation}
\label{eqn:korisno} H_p(F_n)^{(0)} = \big( H_p(SP^{n-1}(M_g),
SP^{n-2}(M_g))\otimes \mathbb{Q}^k \big) \oplus \,
H_p(SP^{n-2}(M_g)) .
\end{equation}
Moreover, ${\rm Im}(\alpha_p) = {\rm Im}(\alpha_p') =
H_p(SP^{n-1}(M_g))\subset H_p(SP^{n}(M_g))$.

By the analysis preceding Proposition~\ref{prop:binom1} we know
that ${\rm Im}(\beta_p)$ is spanned by classes $E_{(I,0)} = e_I =
e_{i_1}\ldots e_{i_p}$ where $1\leq i_1<\ldots < i_p\leq 2g$ and
$p\leq n$. A conclusion is that
\[
{\rm Im}(\alpha_p)\cap {\rm Im}(\beta_p) = {\rm span}\{e_I \mid I
: 1\leq i_1<\ldots < i_p\leq 2g, \, p\leq n-1\}
\]
and a simple calculation completes the proof of the proposition.
\hfill $\square$

\end{document}